\documentclass[11pt,english,leqno]{article}
\usepackage{theorem}
\usepackage{babel}
\usepackage{a4wide}
\usepackage{graphics}
\usepackage{epsfig}
\usepackage{amssymb}
\usepackage{latexsym}
\usepackage{amsmath}

\newtheorem{thh}{Theorem}[section]

\newtheorem{df}[thh]{Definition}

\newtheorem{lem}[thh]{Lemma}
\newtheorem{cor}[thh]{Corollary}
\newtheorem{pro}[thh]{Property}
\newtheorem{prop}[thh]{Proposition}
\newtheorem{rem}[thh]{Remark}

\title{Relations between the leading terms \\ of a polynomial automorphism}
\author{Philippe Bonnet and St\'ephane V\'en\'ereau \\
Mathematisches Institut, Universit\"at Basel \\
Rheinsprung 21, CH-4051 Basel, Switzerland \\
E-mails: philippe.bonnet@unibas.ch, \\
stephane.venereau@unibas.ch}
\date{}

\newcommand{\der}{ \partial}
\newcommand{\dem}{{\em Proof: }}
\newcommand{\qed}{\begin{flushright} $\blacksquare$\end{flushright}}

\newcommand{\KX}{K[x_1,\cdots ,x_n]}
\newcommand{\jj}{\,{\rm j}}

\newcommand{\Ker}{\,{\rm Ker}}

\newcommand{\Z}{\mathbb{Z}}
\newcommand{\N}{\mathbb{N}}
\newcommand{\C}{\mathbb{C}}

\newcommand{\Supp}{\,{\rm Supp}}
\newcommand{\Spec}{\,{\rm Spec}}

\begin{document}
\maketitle

\begin{abstract}
Let $I$ be the ideal of relations between the leading terms of the
polynomials defining an automorphism of $K^n$.  In this paper, we
prove the existence of a locally nilpotent derivation which
preserves $I$. Moreover, if $I$ is principal, i.e. $I=(R)$, we
compute an upper bound for $\deg_2(R)$ for some  degree function
$\deg_2$ defined by the automorphism. As applications, we determine
all the principal ideals of relations for automorphisms of $K^3$ and
deduce two elementary proofs of the Jung-van der Kulk Theorem about the tameness
of automorphisms of $K^{2}$.
\end{abstract}

\section{Introduction}

Throughout this paper, $K$ denotes an algebraically closed field of
characteristic zero. For any positive real numbers $w_1,\cdots, w_n$,
$\deg_1$ stands for the weighted homogeneous degree on the ring
$K[x_1,\cdots,x_n]$ of polynomials in $n$ variables,
which assigns the weights $w_i$ to every $x_i$. We
say that $\deg_1$ defines a positive weighted
homogeneous degree (in short "p.w.h. ").
When
all the weights are equal to $1$, we speak of the (standard)
homogeneous degree on $K[x_1,\cdots,x_n]$.
Given any polynomial $f=\sum a_{k_1,\cdots,k_n} x_1 ^{k_1}\cdots x_n ^{k_n}$ of
degree $d$ for $\deg_1$, its leading term $\overline{f}$ (with respect to $\deg_1$)
is defined as:
$$
\overline{f}=\sum_{k_1 w_1 + \cdots +k_n w_n=d} a_{k_1,\cdots,k_n} x_1 ^{k_1}\cdots x_n^{k_n}
$$
In order to understand polynomial
automorphisms, it is natural to
study the leading terms of their components. More precisely, let
$\Phi=(f_1,\cdots,f_n)$ be a polynomial automorphism of $K^n$.
What can be said about $\overline{f_1},\cdots,\overline{f_n}$?
Except for $n\leq 2$ (see section \ref{tamenesstheorem} below), nothing is known. In this
paper, we would like to determine the algebraic relations
between these leading terms.

\begin{df}
Given a p.w.h.  degree on $K[x_1,\cdots, x_n]$ and an automorphism $(f_1,\cdots,f_n)$
of $K^n$, the set $I$ of polynomials $P\in K[x_1,\cdots,x_n]$ such that $P(\overline{f_1}, \cdots ,\overline{f_n})=0$
is a prime ideal called the {\em ideal of relations}.
\end{df}
Given a p.w.h.  degree $\deg_1$ on $K[x_1,\cdots, x_n]$ and an
automorphism $(f_1,\cdots,f_n)$ of $K^n$, we introduce a new p.w.h.
degree $\deg_2$ on $K[x_1,\cdots,x_n]$, which assigns the weight
$d_i=\deg_1(f_i)$ to every variable $x_i$. Recall that a
$K$-derivation $\der$ of a  $K$-algebra $A$ is locally nilpotent if,
for any $f\in A$, there exists an order $k\geq 0$ with $\der ^k
(f)=0$. Such a derivation is called an {\em LND}. After some simple
observations in section \ref{facts}, we will  establish, in sections
\ref{exLND} and \ref{chute},  the following two properties for
ideals of relations:
\begin{thh} \label{relations}
Let $\Phi=(f_1,\cdots,f_n)$ be an automorphism of $K^n$, let $\deg_1$ be a
p.w.h.  degree on $K[x_1,\cdots,x_n]$ and let $I$ be the ideal of relations.
Then, there exists a nonzero locally nilpotent derivation $\der$ of
$K[x_1,\cdots ,x_n]$ such that
$\der(I)\subset I$. Moreover, if this ideal is principal, i.e. $I=(R)$,
then $\der(R)=0$.
\end{thh}

\begin{thh} \label{inegalite}
Let $\Phi=(f_1,\cdots,f_n)$ be an automorphism of $K^n$. Assume that $\deg_1$ is the
standard homogeneous degree on $K[x_1,\cdots,x_n]$. If the ideal of relations
is principal, i.e. $I=(R)$, then $\deg_2(R)\leq d_1 + \cdots + d_n -n +1$.
\end{thh}
Recall that an automorphism $\tau:K^n \rightarrow K^n$ is elementary if,
up to a permutation of $x_1,\cdots ,x_n$, it can be written as
$\tau(x_1,\cdots,x_n)=(x_1,\cdots, x_n + P)$, where $P\in K[x_1,\cdots,x_{n-1}]$.
 An automorphism of $K^n$ is {\em tame} 
if it is a composition of affine and elementary automorphisms. The Jung-van der Kulk Theorem
states that {\em every automorphism of $K^2$ is tame} (see \cite{Jun,vdK}).
In section \ref{tamenesstheorem}, we give two elementary proofs of this
result, based on Theorems \ref{relations} and \ref{inegalite}.\\
 However
the situation for $n\geq 3$ becomes more complicated. Indeed, there
exist nontame automorphisms in dimension $3$ (see \cite{SU2}). In sections \ref{ssLND}, \ref{bidule}
and \ref{idpr}, applying Theorems \ref{relations} and \ref{inegalite}, we
can determine all possible principal ideals of relations when $n=3$, thereby proving:

\begin{thh}\label{fleche}
Let $\Phi=(f_1,f_2,f_3)$ be an automorphism of $K^3$, and let
$\deg_1$ be the standard homogeneous degree on $K[x_1,x_2,x_3]$. Up
to a permutation of the $f_i$, we may assume that $d_1\leq d_2\leq
d_3$. Assume that the ideal $I$ of relations is principal, i.e.
$I=(R)$. Then there exists a  $\deg_2$-homogeneous
polynomial $h(x_1,x_2)$ whose $\deg_2$ is equal to
$\deg_2(x_3)=\deg_1(f_3)$ or $-\infty$ ($h$ might be zero) such that
$R$ is proportional to one of the following polynomials, where
$x_3':=x_3+h(x_1,x_2)$:
\begin{eqnarray}
&0&\mbox{ (affine case)};\\
 & x_3'  &\mbox{ ("elementarily reducible" case)};\label{t1}\\
 & x_1^{e_1}+cx_2^{e_2}  &\mbox{ with }c\in K^*,\; \gcd(e_1,e_2)=1\label{t2}\;;\\
 & (x_2+ax_1^{e_1})x_3'+cx_1^k   &\mbox{ with }c\in K^*,\,a\in K,\, k\geq 2,\; e_1\geq 1\label{t3}\; ;\\
 & x_1^kx_3'+P(x_1,x_2)   &\mbox{ with }k\geq 1, P\mbox{ $\deg_2$-homogeneous } \label{t4}\; ;\\
  &{x_3'}^2+cx_2^3    &\mbox{ with }c\in K^*\label{t5}\; ;\\
 & {x_3'}^2+cx_1x_2^2  &\mbox{ with }c\in K^*\label{t6}\; ;\\
 & {x_3'}^2+cx_1^{r_1}  &\mbox{ with }c\in K^*,r_1\in 3+2\N\label{t7}\; ;\\
& {x_3'}^2+cx_1^{r_1}x_2  &\mbox{ with }c\in K^*, r_1\geq 1\label{t8}\; ;\\
& {x_3'}^2+ax_1^{e_1}+bx_2^2  &\mbox{ with }ab\in K^*,\; e_1\in 3+2\N\label{t9}\; ;\\
& {x_3'}^2+(ax_1^{e_1}+bx_2)x_1^{r_1} & \mbox{ with }ab\in K^*,\; e_1\geq 1  \label{t8.5} \; ;   \\
& {x_3'}^2+(a_1x_1^{e_1}+b_1x_2)(a_2x_1^{e_1}+b_2x_2) &\mbox{ with }a_1b_2-b_1a_2\neq 0,\;e_1\geq 1\label{t3.5}\; ;\\
 & {x_3'}^2+(a_1x_1+b_1x_2)(a_2x_1+b_2x_2)^2   & \mbox{ with }(a_i, b_i) \neq (0,0) \mbox{ for } i=1,2\label{t6.5}\; ; \\
 & {x_3'}^2+ (ax_1^{e_1}+bx_2)^2 x_1 & \mbox{ with } ab\neq 0,\; e_1\geq 2 \label{t14}\; .
\end{eqnarray}
\end{thh}
Two polynomials $f$ and $g$ of $K[x_1,\cdots ,x_n]$ are said {\em tamely equivalent}
if there exists a tame automorphism $\psi$ of $K^n$ such that $f=g\circ \psi$. With $x_3'$ as in the theorem,
 $R(x_1,x_2,x_3)$ and $R(x_1,x_2,x_3')$ are of course always tamely equivalent.\\
As a consequence of Theorem \ref{fleche}, we will prove in section \ref{democor} the following

\begin{cor}\label{THEcor}
Let $\Phi$ be an automorphism of $K^3$ and let $\deg_1$ be the
standard homogeneous degree on $K[x_1,x_2,x_3]$. If the ideal $I$
of relations is principal, i.e. $I=(R)$, then $R$ is tamely equivalent to
one of the following polynomials: $0$, $x_3$, $x_1 ^r + x_2 ^s$ with $\gcd(r,s)=1$,
$x_1^k x_3  +P(x_1,x_2)$ with $k\geq 1$.
\end{cor}

\section{Some simple but useful facts}\label{facts}
A derivation $\der$ on a general commutative ring $A$ with the same local nilpotence condition as in the introduction
is also called a LND. One may define a function $\deg_\der$ by $\deg_\der(a)=\max\{n|\der^n(a)\neq 0\}$ if $a\neq 0$
 and $\deg_\der(0)=-\infty$. The following lemma is well-known, see e.g. \cite{ML1}:
\begin{lem}\label{degder}
  If $A$ is a domain then $\deg_\der$ is a degree function, i.e. $\deg_\der(a+b)\leq\max\{\deg_\der(a),\deg_\der(b)\}$
  and $\deg_\der(ab)=\deg_\der(a)+\deg_\der(b)$.
 \end{lem}
\begin{cor}\label{da=ba}
  If $A$ is a domain, $a\in A$ and $\der(a)=ba$ for some $b\in A$ then $\der(a)=ba=0$.
\end{cor}
\dem
  If $\der(a)=ba\neq 0$ then $\deg_\der(a)=\deg_\der(\der a)+1=\deg_\der(ab)+1=\deg_\der(a)+\deg_\der(b)+1$
   and one gets $\deg_\der(b)=-1$ which is impossible.
\qed

We will also need to consider the {\em jacobian} (determinant) of $n$ polynomials $P_1,\cdots,P_n\in \KX$:

$$
\jj(P_1,\cdots,P_n)=\det (\partial P_i/\partial x_{j})_{1\leq i,j\leq n}\; .
$$
It is clearly a $K$-derivation in every one of its entries i.e.
every $P_i\mapsto \jj(P_1,\cdots,P_n)$ is  a derivation (when all
the other $P_j$'s are fixed). From its definition, it also clear
that $\deg_1(\jj(P_1,\cdots,P_n))\leq
\deg_1(P_1)+\cdots+\deg_1(P_n)-(\deg_1(x_1)+\cdots+\deg_1(x_n))$. In
parti\-cu\-lar, if $\deg_1$ is the standard homogeneous degree then
\begin{eqnarray}\label{degj<}
\deg_1(\jj(P_1,\cdots,P_n))\leq \deg_1(P_1)+\cdots+\deg_1(P_n)-n\; .
\end{eqnarray}
Note also that, if $\Psi=(h_1,\cdots,h_n)$ is an automorphism, then
$\jj(h_1,\cdots,h_n)=\mu\in K^*$  (well-known fact, consequence of
the chain rule applied on its composition with its inverse) and one
has, $\forall i=1,\cdots,n$,
\begin{eqnarray}
\jj(h_1,\cdots,h_{i-1},P\circ\Psi,h_{i+1},\cdots,h_n) & = &
                                          \jj(h_1,\cdots,h_{i-1},P(h_1,\cdots,h_n),h_{i+1},\cdots,h_n)\nonumber\\
                                                      & = & \jj(h_1,\cdots,h_n)
                                                         \frac{\partial P}{\partial x_i}(h_1,\cdots,h_n)\nonumber\\
\jj(h_1,\cdots,h_{i-1},P\circ\Psi,h_{i+1},\cdots,h_n) & = & \mu\frac{\partial P}{\partial x_i}\circ\Psi\; .\label{jrondpsi}
\end{eqnarray}

Now let $\deg_1$, $\Phi=(f_1,\cdots,f_n)$, $I$ and $\deg_2$ be
as in the introduction. Recall that $\bar f$  denotes the leading
term of the polynomial $f\in K[x_1,\cdots,x_n]$ with respect to
$\deg_1$. Likewise, we will denote $\tilde P$ the leading term of a
polynomial $P\in K[x_1,\cdots,x_n]$ with respect to $\deg_2$.
\begin{lem}\label{1<2}
  Let $P$ be any polynomial in $K[x_1,\cdots,x_n]$. Then $\deg_1(P\circ \Phi)\leq \deg_2 (P)$ and the inequality
  is strict if and only if $\tilde P\in I$ and $P\neq 0$.
\end{lem}
\dem
If $P=0$ this is trivial so let's us assume this is not the case. \\
If $p_{i_1\cdots i_n}$ denotes the coefficient of $x_1^{i_1}\cdots
x_n^{i_n}$ in $P$ then, by definition of the  degree $\deg_2$, there
exists an index $i_1\cdots i_n$ such that $i_1d_1+\cdots
i_nd_n=\deg_2(P)$ and $p_{i_1\cdots i_n}\neq 0$ and one may write:
$$
P(x_1,\cdots ,x_n)  = \sum_{i_1 d_1 + \cdots + i_n d_n=\deg_2 (P)} p_{i_1\cdots i_n} x_1 ^{i_1} \cdots x_n ^{i_n} +
\sum_{i_1 d_1 + \cdots + i_n d_n<\deg_2 (P)} p_{i_1\cdots i_n} x_1 ^{i_1} \cdots x_n ^{i_n}
$$
and, by definition of $\tilde P$,
$$
P(x_1,\cdots,x_n)  =  \tilde P(x_1,\cdots,x_n) + \sum_{i_1 d_1 + \cdots + i_n d_n<\deg_2 (P)}
p_{i_1\cdots i_n} x_1 ^{i_1} \cdots x_n ^{i_n}
$$
Composing with $\Phi=(f_1,\cdots,f_n)$ we get

  \begin{eqnarray*}
    P\circ \phi=P(f_1,\cdots,f_n) & = & \tilde P(f_1,\cdots,f_n) + \sum_{i_1 d_1 + \cdots + i_n d_n<\deg_2 (P)}
    p_{i_1\cdots i_n} f_1 ^{i_1} \cdots f_n ^{i_n}\\
                                   & = & \tilde P(\bar f_1,\cdots,\bar f_n)+ldt\\
  \end{eqnarray*}
  with
$$
  ldt=\tilde P(f_1,\cdots,f_n)-\tilde P(\bar f_1,\cdots,\bar f_n)+\sum_{i_1 d_1 + \cdots + i_n d_n<\deg_2 (P)}
  P_{i_1\cdots i_n} f_1 ^{i_1} \cdots f_n ^{i_n}.
$$ One has
$$
\deg_1(ldt)\leq \max\{\deg_1(\tilde P(f_1,\cdots,f_n)-\tilde P(\bar f_1,\cdots,\bar f_n)),\;
\deg_1(\sum_{i_1 d_1 + \cdots + i_n d_n<\deg_2 (P)} P_{i_1\cdots i_n} f_1 ^{i_1} \cdots f_n ^{i_n})\}
$$ hence $\deg_1(ldt)<\deg_2(P)$. Moreover either $\deg_1(\tilde P(\bar f_1,\cdots,\bar f_n))=\deg_2(P)$ or
 $\tilde P(\bar f_1,\cdots,\bar f_n)=0$.\\
Now let's sum up:
$$
\left\{
\begin{array}{l}
P\circ \Phi=\tilde P(\bar f_1,\cdots,\bar f_n)+ldt\\
\deg_1(\tilde P(\bar f_1,\cdots,\bar f_n))=\deg_2(P)>\deg_1(ldt)\mbox{ or }\tilde P(\bar f_1,\cdots,\bar f_n)=0
\end{array}
\right.
$$
and the result follows from the following
  trivial remark (which holds for any degree function):
  \begin{eqnarray*}
    \deg a+b\leq\max\{\deg a,\deg b\} & \mbox{ and }& \deg a<\deg
    b\Rightarrow \deg a+b=\deg b\;
  \end{eqnarray*}
applied to $\deg=\deg_1$, $a=\tilde P(\bar f_1,\cdots,\bar f_n)$ and $b=ldt$.
\qed
Note that, in the case when $\deg_1$ is the standard homogeneous degree, the weights of $\deg_2$, namely the
$d_i=\deg_1(f_i)$'s are all at least $1$, since $\Phi=(f_1,\cdots,f_n)$ is an automorphism. Consequently one gets the
 implication: $\deg_1(P)>1$ $\Rightarrow$ $\deg_2(P)>1$.
\begin{cor}\label{degPPhi=1} If $\deg_1$ is the standard homogeneous degree and $P\in K[x_1,\cdots,x_n]$  is such that
 $\deg_1(P\circ \Phi)=1$ then $\tilde P\in I$ unless $P$ is affine i.e. $\deg_1(P)\leq 1$.
\end{cor}
\dem
If $P$ is not affine then $\deg_2(P)>1=\deg_1(P\circ \Phi)$ and it follows from Lemma \ref{1<2} that $\tilde P\in I$.
\qed
\begin{cor}\label{I=0} If $\deg_1$ is the standard homogeneous degree
  then $I=(0)$ if and only if $\Phi$ is an affine automorphism i.e. $\deg_1(f_i)=1$, $\forall i=1,\cdots,n$.
\end{cor}
\dem
  If $\Phi$ is affine then
$$
\bar \Phi:=(\bar f_1,\cdots,\bar f_n)=(f_1-f_1(0,\cdots,0),\cdots,f_n-f_n(0,\cdots,0))=\tau\circ\Phi
$$
 where $\tau$ is the automorphism (translation) $(x_1-f_1(0,\cdots,0),\cdots,f_n-f_n(0,\cdots,0))$ also $\bar\Phi$
  it self is an automophism preventing that $P\circ\bar\Phi=P(\bar f_1,\cdots,\bar f_n)=0$ for any $P\neq 0$. It follows
  that $I=(0)$.\\
Assume now that $I=(0)$.  Let $\Phi^{-1}=(g_1,\cdots,g_n)$ be the
inverse of $\Phi$. One has $\forall i=1,\cdots,n$, $g_i\circ
\Phi=x_i$ hence $1=\deg_1(x_i)=\deg_1(g_i\circ \Phi)$ and $0\neq
\tilde g_i\notin I$. By Corollary \ref{degPPhi=1}, $g_i$ is affine
$\forall i=1,\cdots,n$, hence so are $\Phi^{-1}=(g_1,\cdots,g_n)$
and $\Phi$. \qed Remark that for p.w.h. degrees other than the
standard one, none of the implications in Corollary \ref{I=0} holds.
Indeed, take the weights $w_1=\deg_1(x_1)=1,w_2=\deg_1(x_2)=3$. Then
$(f_1,f_2)=(x_1+x_2,x_2)$ is an affine  automorphism but
$\bar f_1=x_2$ and $\bar f_2=x_2$ fullfill the non trivial relation $\bar f_1-\bar f_2=0$.\\
Conversely $(f_1,f_2)=(x_1,x_2+x_1^2)$ is a non-affine automorphism but $\bar f_1=x_1$ and $\bar f_2=x_2$ are
 algebraically independant.\\

Finally, we point out an easy fact linking degree and partial
derivative (here $\frac{\partial}{\partial x_n}$) in  the following
\begin{lem}\label{deg2der}
  Let $P$ be a polynomial in $\KX$. If $\frac{\partial\tilde P}{\partial x_n}\neq 0$ then
  $\widetilde{\frac{\partial  P}{\partial x_n}}=\frac{\partial \tilde P}{\partial x_n}$ and
   $\deg_2(\frac{\partial  P}{\partial x_n})=\deg_2(P)-d_n$. More generally,
    while $\frac{\partial^k\tilde P}{\partial x_n^k}\neq 0$ one has
    $\widetilde{\frac{\partial^k  P}{\partial x_n^k}}=\frac{\partial^k \tilde P}{\partial x_n^k}$ and
    $\deg_2(\frac{\partial^k  P}{\partial x_n^k})=\deg_2(P)-kd_n$.
\end{lem}
\dem
  Derivating  $P=\tilde P+ldt$ where $\deg_2(ldt)<\deg_2(P)=\deg_2(\tilde P)$ one gets
\begin{eqnarray*}
\frac{\partial P}{\partial x_n} & = & \frac{\partial \tilde P}{\partial x_n}+\frac{\partial ldt}{\partial x_n}\; .
\end{eqnarray*}
 Since $\tilde P$ is $\deg_2$-homogeneous, so is  $ \frac{\partial \tilde P}{\partial x_n}$ and its  $\deg_2$ is equal
 to $\deg_2(P)-\deg_2(x_n)=\deg_2(P)-d_n$. On the other hand, one has $\deg_2(\frac{\partial ldt}{\partial x_n})\leq
 \deg_2(ldt)-d_n<\deg_2(P)-d_n$ and everything comes from the equality above.
\qed


\section{Proof of the first Theorem}\label{exLND}

Let $\der$ be a $K$-derivation of $K[x_1,\cdots ,x_n]$. Given a p.w.h.  degree
$\deg$ on $K[x_1,\cdots ,x_n]$, the degree of $\der$ is defined as:
$$
\deg(\der)=\sup\{\deg(\der(P)) - \deg(P)|\; P\in
K[x_1,\cdots ,x_n]-\{0\}\}
$$
Note that this degree is upper-bounded, i.e. $\deg(\der)<+\infty$, and equals
$-\infty$ if and only if $\der$ is the zero derivation. To see this,
it suffices to check that:
$$ \deg(\der)=\sup\{\deg(\der(x_i)) - \deg(x_i)|\; i=1,\cdots ,n\}$$
By construction, we have the
following inequality:
$$
\forall P \in K[x_1,\cdots ,x_n], \quad \deg(\der(P))\leq \deg(\der) +
\deg(P)\; .
$$
The leading term $\overline{\der}$ of a $k$-derivation $\der$ is the
$k$-linear endomorphism of $K[x_1,\cdots ,x_n]$, which maps every
weighted homogeneous polynomial $P$ to the weighted homogeneous term
of $\der(P)$ of degree $\deg(\der) + \deg(P)$. By construction,
$\overline{\der}$ is a $k$-derivation of $K[x_1,\cdots ,x_n]$.
More precisely, assume that $\der= \sum_i a_i \der/\der x_i$, where
every $a_i$ is a polynomial. If $r$ is the degree of $\der$ and if
$b_i$ denotes the weighted homogeneous term of degree $r + \deg(x_i)$ of $a_i$,
then:
$$
\overline{\der} = \sum_{i=1} ^n b_i \frac{\der}{\der x_i}\; .
$$
Let $\Phi=(f_1,\cdots ,f_n)$ be a polynomial automorphism of $K^n$, with
jacobian equal to $\lambda\in k$, and denote by $\Phi^{-1}=(g_1,\cdots ,g_n)$
its inverse (also with jacobian $\lambda^{-1}$). Let $\deg_1$ be a p.w.h.  degree, which assigns
the weight $w_i$ to each variable $x_i$, and let $\deg_2$ be the
degree defined in the introduction. For any $i=1,\cdots ,n$, consider
the $k$-derivations $\delta_i$ and $\Delta_i$ on $\KX$, defined
for any polynomial $P$ by the formulas:
$$
\delta_i(P)=\lambda^{-1} \frac{\der P}{\der x_i} \quad \mbox{and} \quad
\Delta_i(P)=\jj(g_1,\cdots ,g_{i-1},P,g_{i+1},\cdots ,g_{n})\; .
$$
Note that each derivation $\delta_i$ has degree $-w_i$ with respect to
$\deg_1$. Theorem \ref{relations} will be a straightforward consequence of the following
lemmas, where we prove that the leading term of one of the $\Delta_i$ for
$\deg_2$ is an LND which stabilizes the ideal $I$ of relations.

\begin{lem} \label{LND01}
For any index $i$ of $\{1,\cdots ,n\}$ and any polynomial $P$ of $\KX$, we have $\Delta_i(P\circ
\Phi^{-1})=\delta_i(P)\circ \Phi^{-1}$ and $\Delta_i(P)\circ
\Phi=\delta_i(P\circ \Phi)$. Moreover, every derivation $\Delta_i$ is locally
nilpotent.
\end{lem}
\dem For any index $i=1,\cdots ,n$ and any polynomial $P$ of $\KX$ we have, by equality (\ref{jrondpsi}):
$$
\Delta_i(P\circ\Phi^{-1})  =  \jj (g_1,\cdots ,g_{i-1},P\circ \Phi^{-1},g_{i+1},\cdots ,g_{n})=\lambda^{-1}\frac{\der P}{\der x_i}\circ \Phi^{-1}=\delta_i(P)\circ \Phi^{-1}\; .
$$
Now, for any polynomial $P$, if we set $Q=P\circ \Phi$, then we find:
$$
\Delta_i(P)\circ \Phi=\Delta_i(Q\circ \Phi^{-1})\circ \Phi=\delta_i(Q)\circ
\Phi^{-1}\circ \Phi = \delta_i(Q)=\delta_i(P\circ \Phi)\; .
$$
Using the second formula, one can prove by induction on $k\geq 1$ that
$\Delta_i ^k (P)\circ \Phi=\delta_i ^k (P\circ \Phi)$ for any $P$. Now fix a
polynomial $P$. Since $\delta_i$ is locally nilpotent, there exists an
order $r\geq 1$ such that $\delta_i ^r(P\circ \Phi)=0$. In particular
$\Delta_i ^r(P)\circ \Phi=0$. Since $\Phi$ is an automorphism, $\Delta_i ^r(P)=0$.
Since this holds for any $P$, $\Delta_i$ is locally
nilpotent. \qed

\begin{lem} \label{ineg}
There exists an index $i$ such that $\deg_2(\Delta_i)\geq \deg_1(\delta_i)=-w_i$.
\end{lem}
\dem Suppose on the contrary that $\deg_2(\Delta_i)<\deg_1(\delta_i)=-w_i$
for any $i=1,\cdots ,n$. Fix an index $j$ for which $d_j=\min_i\{d_i\}$. Using Lemma
\ref{LND01} and the fact that $x_j=f_j\circ \Phi^{-1}$, we obtain for any index $i$:
$$
\deg_2(\Delta_i(x_j))=\deg_2(\Delta_i(f_j\circ \Phi^{-1}))=\deg_2(\delta_i(f_j)\circ \Phi^{-1})<d_j -w_i
$$
Since $d_j$ is minimal among all the $d_i$, every non-constant polynomial has
degree $\geq d_j$ with respect to $\deg_2$. Since every weight $w_i$ is positive,
$d_j -w_i<d_j$ for any $i$. So $\delta_i(f_j)\circ \Phi^{-1}$ is constant for any $i$.
Since $\Phi$ is an automorphism, $\delta_i(f_j)$ is constant for any $i$. This implies
there exists some constants
$a,a_1,\cdots ,a_n$ such that:
$$
f_j=a + a_1 x_1 + \cdots  + a_nx_n
$$
In particular, $d_j$ is the maximum of the $w_i$ for which $a_i\not=0$. Now fix an index
$l$ such that $d_j=w_l$ and $a_l\not=0$. As before, we obtain:
$$
\deg_2(\Delta_l(x_j))=\deg_2(\delta_l(f_j)\circ \Phi^{-1})<d_j -w_l=0
$$
Since every weight $w_i$ is positive, the polynomial $\delta_l(f_j)\circ \Phi^{-1}$ must
be equal to zero. But this is impossible since $\delta_l(f_j)=a_l\not=0$.
\qed

\begin{lem}
Let $i$ be an index such that $\deg_2(\Delta_i)\geq -w_i$. Then the leading part $\overline{\Delta_i}$
of $\Delta_i$ with respect to $\deg_2$ is locally nilpotent and stabilizes the ideal $I$.
\end{lem}
\dem Let $i$ be an index such that $\deg_2(\Delta_i)\geq
-w_i=\deg_1(\delta_i)$. Since $\Delta_i$ is locally nilpotent  by
Lemma \ref{LND1}, its leading term $\overline{\Delta_i}$ is also
locally nilpotent (see \cite{vdE}). So we only need to prove that
$\overline{\Delta_i}$ stabilizes $I$. Since the $\overline{f_i}$ are
$\deg_1-$homogeneous of degrees $d_1,\cdots ,d_n$ , $I$ is a
weighted homogeneous ideal for $\deg_2$. Let $P$ be any nonzero
$\deg_2$-homogeneous element of $I$. Using successively Lemma
\ref{1<2}, the definition of the degree of a derivation and Lemma
\ref{LND01}, we find:
$$
\begin{array}{ccl}
\deg_2(\Delta_i)+ \deg_2(P) & > & \deg_1(\delta_i)+ \deg_1(P\circ \Phi) \\
& & \\
& > & \deg_1(\delta_i(P\circ \Phi)) \\
& & \\
& > & \deg_1(\Delta_i(P)\circ \Phi)
\end{array}
$$
If $\deg_2(\Delta_i(P))< \deg_2(\Delta_i)+ \deg_2(P)$, then $\overline{\Delta_i}(P)=0$
and there is nothing to prove. So assume that $\deg_2(\Delta_i(P))= \deg_2(\Delta_i)+ \deg_2(P)$.
By the inequality proved above, we get:
$$
\deg_2(\Delta_i(P))> \deg_1(\Delta_i(P)\circ \Phi)
$$
Therefore, $\overline{\Delta_i(P)}$ belongs to $I$ by Lemma
\ref{1<2}. Since $\deg_2(\Delta_i(P))= \deg_2(\Delta_i)+
\deg_2(P)$, we have $\overline{\Delta_i(P)}=\overline{\Delta_i}(P)$
by definition of the leading term of a derivation. So
$\overline{\Delta_i}(P)$ belongs to $I$. Since this holds for any
element $P$ of $I$, the derivation $\overline{\Delta_i}$ stabilizes
$I$. \qed The last part of Theorem \ref{relations} ("Moreover, if
this ideal ...") follows from Corollary \ref{da=ba}.


\section{On the parachute for the degree} \label{chute}

In this section, we prove Theorem \ref{inegalite} using the so-called  "parachute" defined in  \cite{V07}. \\
Since one wants to majorate $\deg_2(R)$, one may assume that $R\neq
0$. Moreover $R$ cannot be constant so at least one  of its partial
derivatives $\frac{\partial R}{\partial x_i}$ is nonzero. Without
loss of generality we will assume that  $\frac{\partial R}{\partial
x_n}\neq 0$. Finally recall that $I=(R)$ must be prime so $R$ is
irreducible.
 Let us list our assumptions:
\begin{itemize}
\item $\deg_1$
denotes the standard homogeneous degree on $K[x_1,\cdots,x_n]$;
\item $\Phi=(f_1,\cdots,f_n)$ is an automorphism of $K[x_1,\cdots,x_n]$;
\item $\deg_2$ is the p.w.h. degree defined by $\deg_2(x_i)=d_i=\deg_1(f_i)$, $\forall i=1,\cdots,n$;
\item $\bar P$ resp. $\tilde P$ denotes the leading term of a polynomial with respect to $\deg_1$ resp. $\deg_2$;
\item the ideal of relations  is principal i.e. $I=(R)$, for some irreducible $R\in K[x_1,\cdots,x_n]$ such that   $
\frac{\partial R}{\partial x_n}\neq 0$.
\end{itemize}

  In our situation the {\em parachute} of $f_1,\cdots,f_n$, denoted $\nabla$, is simply defined as follows:
$$
  \nabla:=d_1+\cdots+d_n-n\;\;.
$$
The parachute prevents the degree of a polynomial from "falling too much":
\begin{pro} \label{fallschirm}
For any polynomial $P$  in
$K[y_1,\cdots,y_m]$, one has the following minoration:
$$
\deg_1 (P \circ \Phi) \geq \deg_1(\frac{\partial P}{\partial x_n}\circ \Phi) +d_n-\nabla\mbox{ and, inductively,}
$$
$$
\deg_1 (P \circ \Phi) \geq \deg_1(\frac{\partial^k P}{\partial x_n^k}\circ \Phi) +kd_n-k\nabla, \;\forall k\geq 0\;.
$$

\end{pro}
\dem  It follows from equality (\ref{jrondpsi}) applied to
$\Phi=(f_1,\cdots,f_n)$ that $\deg_1(\jj(f_1,\cdots,f_{n-1},
P\circ\Phi))=\deg_1 (\frac{\partial P}{\partial x_n}\circ\Phi)$. On
the other hand, by inequality (\ref{degj<}), one has
$\deg_1(\jj(f_1,\cdots,f_{n-1},P\circ\Phi))\leq
\deg_1(f_1)+\cdots+\deg_1(f_{n-1})+\deg_1(P\circ\Phi)-n$ so we have
\begin{eqnarray*}
  \deg_1 (\frac{\partial P}{\partial x_n}\circ\Phi) & \leq & d_1+\cdots+d_{n-1}+\deg_1(P\circ\Phi)-n\\
\deg_1 (\frac{\partial P}{\partial x_n}\circ\Phi)-(d_1+\cdots+d_{n}-n)+d_n & \leq & \deg_1(P\circ\Phi) \\
\deg_1 (\frac{\partial P}{\partial x_n}\circ\Phi)+d_n-\nabla & \leq & \deg_1(P\circ\Phi)
\end{eqnarray*}
as required.
\qed

\begin{lem} \label{ordre}
Let $H$ be a nonzero polynomial in $\KX$ and $k\in\N$ be such that $H\in (R^k)\setminus(R^{k+1})$.
Then  $\frac{\der ^k H}{\der x_n ^k}\notin(R)$.
\end{lem}
\dem We proceed by induction on $k$. For $k=0$ there is nothing to prove. \\
Assume this holds for some $k-1\geq 0$. One has $H=SR^k$ for some $S\notin R$ and
$$
  \frac{\partial H}{\partial x_n}=\frac{\partial S}{\partial x_n}R^k+SR^{k-1}\frac{\partial R}{\partial x_n}
$$
is clearly in $(R^{k-1})$ but not in $(R^k)$ because otherwise
$S\frac{\partial R}{\partial x_n}$ would be in $(R)$  which is
impossible for $S$ by assumption, and for $\frac{\partial
R}{\partial x_n}$ for a degree reason so impossible for the product
since $(R)$ is prime. The induction assumption concludes. \qed

\begin{prop}\label{above}
  Let $P$ be a nonzero polynomial in $\KX$ and $k\in\N$ be such that $\tilde P\in (R^k)\setminus(R^{k+1})$. One has
\begin{description}
\item[{\rm (i)}]  $\deg_1(\frac{\partial^k P}{\partial x^k}\circ\Phi)=\deg_2(\frac{\partial^k P}{\partial x_n^k})=
\deg_2(P)-kd_n$;
\item[{\rm (ii)}]$\deg_1(P\circ \Phi)\geq\deg_2(P)-k\nabla$;
\item[{\rm (iii)}]$\deg_1(P\circ \Phi)\geq k(\deg_2(R)-\nabla)$.
 \end{description}
\end{prop}
\dem
 By Lemma \ref{ordre}, $\frac{\partial^k \tilde P}{\partial x^k}\notin (R)$. In particular
$\frac{\partial^k \tilde P}{\partial x^k}\neq 0$ hence, by Lemma
\ref{deg2der}, $\widetilde{\frac{\partial^k  P}{\partial
x^k}}=\frac{\partial^k \tilde P}{\partial x^k}\notin (R)=I$ and
$\deg_2(\frac{\partial^k P}{\partial x_n^k})=\deg_2(P)-kd_n$ also,
by Lemma \ref{1<2}, $\deg_1(\frac{\partial^k  P}{\partial
x^k}\circ\Phi)=\deg_2(\frac{\partial^k  P}{\partial
x^k})=\deg_2(P)-kd_n$, thereby proving (i). Assertion (ii) is a
direct consequence of (i) and Property \ref{fallschirm}. To prove
(iii) observe that $\tilde P\in (R^k)$ implies that
$k\deg_2(R)\leq\deg_2(\tilde P)=\deg_2(P)$. \qed

{\it Proof of Theorem \ref{inegalite}}: Let $\Phi^{-1}=(g_1,\cdots,
g_n)$ be the inverse of $\Phi$. Since $R\not=0$,  it follows from
Corollary \ref{I=0} that the automorphism $\Phi$ is not affine and
neither is its inverse. In particular, there exists an index $j$
such that $g_j$ is not affine. By definition, $x_j=g_j\circ\Phi$
hence $\deg_1(g_j\circ\Phi)=1$ and, by Corollary \ref{degPPhi=1},
$\tilde g_j\in I=(R)$. It follows that the integer $k$ such that
$\tilde g_j\in (R^k)\setminus (R^{k+1})$ is positive. Now applying
Proposition \ref{above} (iii) with $P=g_j$ one has
$1=\deg_1(g_j\circ\Phi)\geq k(\deg_2(R)-\nabla)\geq
\deg_2(R)-\nabla$ and consequently $\deg_2(R)\leq
1+\nabla=d_1+\cdots+d_n-n+1$. \qed

\section{On the Jung-van der Kulk Theorem} \label{tamenesstheorem}

In this section, we are going to give two elementary proofs of
the Jung-van der Kulk Theorem. This theorem states that
{\em every automorphism $\Phi$ of $K^2$ is tame, i.e. $\Phi$ is a
composition of affine and elementary automorphisms}. Both proofs
are based on the computation of the ideal of relations for
$\Phi$, and they use induction on the sum $n=\deg_1(f) + \deg_1(g)$,
where $\deg_1$ is the standard homogeneous degree on $K[x,y]$
and where $\Phi=(f,g)$.

They proceed as follows. First, if $n=2$,
then $\Phi$ is affine and there is nothing to show. So assume
that the theorem holds at every order from $2$ to $n$, and let
$\Phi$ be an automorphism such that $\deg_1(f) + \deg_1(g)=n+1$.
Up to a permutation of $f$ and $g$ (which corresponds to composing
$\Phi$ with a linear map), we may assume that $\deg_1(f)\geq \deg_1(g)$.
Since $\Phi$ ist not affine, by Corollary \ref{I=0} the ideal $I$ of
relations is distinct from $(0)$, and its height is either equal to $1$
or $2$. Note that, if $I$ were of height $2$, then it would be equal to
$(x,y)$ because $I$ is prime in $K[x,y]$ and weighted homogeneous.
But then $\overline{f}=\overline{g}=0$, which is impossible. So
$I$ is prime of height $1$, hence it is principal. Write $I=(R)$,
where $R$ is irreducible and weighted homogeneous. Assume now
that \\

$(*)$\centerline {\em $R$ is of the form $x-cy^r$, where $c\in K$ and $r\in \mathbb{N}^*$.}\\

Then, this implies that: $\overline{f} - c\overline{g}^r=0 \quad \mbox{and} \quad \deg_1(f-cg^k)<\deg_1(f)$.
Consider the map $\Phi'= (f',g')$, where $f'=f-cg^k$ and $g'= g$. As a composition
of $\Phi$ with an elementary automorphism, it is itself an automorphism of $K^2$.
Since $\deg_1(f')+\deg_1(g')\leq n$, $\Phi'$ is tame by the induction's hypothesis.
So $\Phi$ is tame and the result follows. Therefore, we only need to prove the
assertion $(*)$. In the following subsections, we exhibit two proofs of $(*)$
as an application of our main Theorems.

\subsection{First proof of the assertion $(*)$}

Let $R$ be an irreducible weighted homogeneous polynomial such that $I=(R)$.
Then $R$ is of the form $R(x,y) = \lambda x ^{s} + \mu y ^{r}$, where $s,r$
are coprime. By Theorem \ref{relations}, there exists a nonzero
weighted homogeneous locally nilpotent derivation $\der$ of $K[x,y]$ such that
$\der(I)\subset I$. In particular, we have $\der(R)=0$. One could directly conclude using \cite{ML2} but we prefer to give a self-contained proof. Write $\der$ as:
$$
\der=a\frac{\der}{\der x} + b \frac{\der}{\der y}
$$
where $a,b$ are polynomials. Let $R^n$ be the highest power of $R$ that divides both
$a$ and $b$. Since $\der$ is locally nilpotent and that $\der(R)=0$, $D=\der /R^n$ is
also locally nilpotent. Since $D(R)=0$ and that $a/R^n, b/R^n$ do not both belong to
$(R)$, $D$ induces a nonzero locally nilpotent derivation $D'$ on the ring:
$$
A= \frac{K[x,y]}{(R)} = \frac{K[x,y]}{(\lambda x^{s} + \mu y ^{r})} \simeq K[t^{r},t^{s}]
$$
Note that, since $\deg_1(f)\geq \deg_1(g)$, we have $r\geq s$. Since
$s,r$ are coprime, the integral closure of  $K[t^{r},t^{s}]$ is
equal to $K[t]$. By Seidenberg's Theorem (see \cite{Sei}), $D'$
extends to a derivation of $K[t]$. By Vasconcelos' Theorem (see
\cite{Vas}), this extension is locally nilpotent. So $D'$ is of the
form $\theta \der/\der t$, where $\theta$ belongs to $K^*$. In
particular, $D'(t^{s})=s\theta t^{s-1}$ belong to $K[t^{r},t^{s}]$.
Since $r\geq s$, this is only possible if $s=1$, and $R$ is of the
form $\lambda x+ \mu y^r$.

\subsection{Second proof of the assertion $(*)$}

Let $R$ be an irreducible weighted homogeneous polynomial such that $I=(R)$.
Then $R$ is of the form $R(x,y) = \lambda x ^{s} + \mu y ^{r}$, where $s,r$
are coprime. Set $\deg_1(f)=d_1$ and $\deg_1(g)=d_2$. By assumption, we have
$d_2 \leq d_1$. By applying Theorem \ref{inegalite} to the polynomial $R$,
we find:
$$
s d_ 2 \leq d_1 + d_2 -2 < 2d_2
$$
Therefore, $s$ is either equal to $0$ or $1$. But $s$ cannot be equal to $0$,
otherwise $\overline{g}^r$ would be equal to zero, which is impossible. Therefore,
$s=1$ and $R(x,y) = \lambda x + \mu y ^{r}$.


\section{A few surfaces without LND}\label{ssLND}
By Theorem \ref{relations}, in order to show that some polynomials
do not appear as generators of relations between  the leading terms
of an automorphism it is sufficient to show that they are not
annihilated by a nontrivial LND. Actually we will show a bit more
about those polynomials, namely that the regular functions ring on
the surfaces they define admits in turn no nontrivial LND, and use
the following
\begin{lem} \label{LND1}
Let $\deg_1$ be a p.w.h.  degree on $K[x_1,\cdots ,x_n]$. Let $R$ be
a $\deg_1-$homogeneous irreducible  element of $K[x_1,\cdots
,x_n]-K$. Assume there exists a nonzero locally nilpotent derivation
$D$ on $K[x_1,\cdots ,x_n]$ such that $D(R)=0$. Then the quotient
ring $A=K[x_1,\cdots ,x_n]/(R)$ admits a nonzero weighted
homogeneous locally nilpotent derivation.
\end{lem}
\dem Since $D$ is locally nilpotent, its leading term $\overline{D}$
is nonzero and locally nilpotent (see \cite{vdE}).  Since $R$ is
weighted homogeneous and that $D(R)=0$, we have $\overline{D}(R)=0$.
So we may assume that $D$ is weighted homogeneous. Write $D$ as:
$$
D = \sum_{i=1} ^n a_i\frac{\der}{\der x_i}\; .
$$
Let $R^k$ be the maximal common power of $R$ dividing $a_1,\cdots
,a_n$. Since $D(R)=0$, the derivation $\der=D/R^k$ is a locally
nilpotent derivation on $K[x_1,\cdots ,x_n]$ which induces a
weighted homogeneous locally nilpotent derivation $\der '$ on the
quotient ring $A$. Moreover, since not all $a_i/R^k$ are divisible
by $R$, $\der '$ is nonzero. \qed

Using some results and technics of \cite{KZ00}  we prove
\begin{prop}\label{noLND}
The quotient algebra $K[x_1,x_2,x_3]/(R)$ admit no nonzero LND when
$R$ is one of the following polynomials:
$$
\begin{array}{ccc}
x_3^2+ax_1^{4}+bx_2^{3} &  \mbox{where} & ab\neq 0\\
x_3^2+ax_1^{5}+bx_2^{3} &  \mbox{where} & ab\neq 0\\
x_3^2+(a_1x_1+b_1x_2)(a_2x_1+b_2x_2)(a_3x_1+b_3x_2) & \mbox{where} & a_ib_j - a_jb_i\not=0, \forall i\not=j \\
 x_3^2+(a_1x_1^{e_1}+b_1x_2)(a_2x_1^{e_1}+b_2x_2)x_1 & \mbox{where} & a_1b_2 - a_2b_1\not=0,\; e_1\geq 0\\

x_3 ^2 + (ax_1 ^3 + bx_2 ^2)x_2 & \mbox{where} & ab\neq 0\; .
\end{array}
$$
\end{prop}
\dem Up to multiplication of $x_1$ and $x_2$ by appropriate scalars
and renaming the indeterminates, our first two polynomials define
the Platonic surfaces $S_{2,3,4}$ and $S_{2,3,5}$ where
$S_{k,l,m}=\{x^k+y^l+z^m=0\}\subset K^3$ as in the proof of Lemma 4
in \cite{KZ00}. Even if, in the cited article, $K$ is the complex
field, the proofs that $S_{2,3,4}$ and $S_{2,3,5}$ admit no
non-trivial LND still works word for word for a field as general as
ours.  Our third polynomial defines a surface isomorphic to
$S_{2,3,3}$: in a word one chooses a linear change of the
coordinates $x_1,x_2$ corresponding to the projective map taking
$(a_i,b_i)$ to $(1,\lambda^i)$ for $i=1,2,3$ where $\lambda$ is a
primitive cubic root of the unity and one rescales $x_3$. Again by
\cite{KZ00}, $S_{2,3,3}$ admit no non-trivial LND.\\
Now we assume that $R$ is one of the two last polynomials with
weights
$(\deg_1(x_1),\deg_1(x_2),\deg_1(x_3))=(d_1,d_2,d_3)=(2,2e_1,2e_1+1)$
resp. $(4,6,9)$, and
that the quotient ring $A=K[x_1,x_2,x_3]/(R)$ admits a nonzero
 LND $\der$. Up to taking the leading term of $\der$ as in
 $\cite{KZ00}$ (see also the beginning of our section \ref{exLND})
 one may assume that this derivation is weighted homogeneous.

First, we are going to prove that $\ker \der$ contains an
irreducible weighted homogeneous polynomial $P(x_1,x_2)$.  Since $A$
has transcendence degree $2$ over $K$, there exists a nonconstant
element $f$ of $A$ such that $\der(f)=0$. Since $A$ is graded and
that $\der$ is weighted homogeneous, we may assume that $f$ is also
weighted homogeneous. Consider $R$ as a polynomial of degree $2$ in
$x_3$. After division by $R$, we may assume that $f$ is of the form:
$$
f=f_1(x_1,x_2)x_3+f_0(x_1,x_2)
$$
where $f_1,f_0$ are polynomials. Since $\deg_1(x_i)$ is even for
$i=1,2$, $\deg_1(f_1)$ and $\deg_1(f_2)$ are even. Since $f$ is
 weighted homogeneous and  $\deg_2(x_3)$ is odd,
then either $f = f_1(x_1,x_2)x_3$ or $f=f_0(x_1,x_2)$. In all cases,
since $x_3 ^2$ is a polynomial in $x_1,x_2$,  $f^2$ belongs to $\ker
\der \cap K[x_1,x_2]$. So $\ker \der$ contains a nonconstant
weighted homogeneous element $P$ of $K[x_1,x_2]$. Since $\der$ is
locally nilpotent, $\ker \der$ is factorially closed (well-known fact, see e.g. \cite{ML1}) and every
irreducible factor of $P$ in $K[x_1,x_2]$ lies in $\ker \der$. So we may assume that $P$ is irreducible.

Second, note that, since $P$ is weighted homogeneous, it is either
equal (up to a scalar multiple) to $x_1$, or $x_2$, or $cx_1 ^r +
x_2 ^s$, where $c\in k$. Now
let us distinguish between the two cases.\\

If $R=x_3^2+(a_1x_1^{e_1}+b_1x_2)(a_2x_1^{e_1}+b_2x_2)x_1 $ with the
weights $(2,2e_1,2e_1+1)$ and $P=cx_1 ^r + x_2 ^s$ then $r=e_1$ and
$s=1$ so that it is weighted homogeneous. The change of coordinates
$(x_1,x_2,x_3)\mapsto (x_1,x_2-cx_1^r,x_3)$ doesn't affect the form
of $R$, so we may assume $c=0$. So we are left with the two cases
$x_1\in \ker \partial$ or $x_2\in \ker \partial$.

 Assume that $\der(x_1)=0$, and let $L$ be
the algebraic closure of the field $K(x_1)$. Define a new derivation $D$
on the ring $B=L \otimes_{K[x_1]} A$ as:
$$
D : L \otimes_{K[x_1]} A \longrightarrow L \otimes_{K[x_1]} A, \quad
a \otimes b \longmapsto a \otimes \der(b)\; .
$$
Since $\der$ is nonzero locally nilpotent, $D$ is also nonzero
locally nilpotent. Moreover, $B$ is isomorphic to the ring
$L[x_2,x_3]/(x_3^2+x_2^{2}-1)$. Since $\Spec(B)$ is isomorphic to
$L^*$, $D$ is equal to zero, which is impossible. So
$\der(x_1)\not=0$ and $\partial (x_2)=0$. Then, following the same
construction, we obtain a nonzero locally nilpotent derivation $D$
on the coordinate ring of the curve $x_3^2+(a_1x_1^{e_1}+b_1x_2)(a_2x_1^{e_1}+b_2x_2)x_1=0 $
 over $L$, where $L$ is the
algebraic closure of the field $K(x_2)$. Since this curve is smooth
and hyperelliptic, its geometric genus is equal to $e_1 - 1$. In
particular, this curve cannot be isomorphic to $L$ unless $e_1=1$.
But in this case we get $L^*$; therefore $D$ is
equal to $0$, hence a contradiction.\\

If $R=x_3 ^2 + (ax_1 ^3 + bx_2 ^2)x_2$ with the weights $(4,6,9)$
then the case $P=cx_1 ^r + x_2 ^s$ is excluded since $(r,s)$ would
have to be $(3,2)$ which is impossible by  \cite{ML2}. It follows
that $x_1$ or $x_2$ is in $\Ker\der$ and we get a contradiction as
in the previous case.

\qed


\section{A special surface without LND}\label{bidule}

In this section, we are going to establish the following result:

\begin{prop}\label{special}
Let $e_1$ be in $3+2\N$, $a,b\in K^*$ and
$R=x_3^2+(ax_1^{e_1}+bx_2^{2})x_1$. Then the quotient algebra
$K[x_1,x_2,x_3]/(R)$ admit no nonzero LND.
\end{prop}
Up to replacing each $x_i$ with $\lambda_i x_i$ for some suitable constants $\lambda_i$,
we may assume that $a=1$ and $b=-1$. For any element $e_1$ of $3+2\N$, denote by $S$ the surface
of $K^3$ given by the equation $R=0$. Consider the Danielewski surface $S'$ of $K^3$,
given by the equation:
$$
\alpha ^{2e_1} - 4\beta \gamma = 0
$$
together with the involution $\sigma$ of $S'$ defined by the formula
$\sigma(\alpha,\beta,\gamma)=(-\alpha,\gamma,\beta)$.  Denote by
$\sigma^*$ the involution of ${\cal{O}}(S')$ defined by
$\sigma^*(f)=f\circ\sigma$. It is easy to check that $\sigma$ has
$(0,0,0)$ as a unique fixpoint, and that the map:
$$
F: S' \longrightarrow S, \quad (\alpha, \beta, \gamma) \longmapsto (\alpha^2, \beta+\gamma, \alpha(\beta - \gamma))
$$
is a well-defined morphism, i.e. $F(S')\subseteq S$. Since $F$ is dominant, it induces
an injective algebra morphism from ${\cal{O}}(S)$ to ${\cal{O}}(S')$. We can therefore
consider ${\cal{O}}(S)$ as a subring of ${\cal{O}}(S')$. Note that $S'$ is weighted
homogeneous for the weights $1, e_1, e_1$. Moreover, the coordinate functions of $F$
are all weighted homogeneous of degree $2, e_1, e_1 + 1$. We begin with a few lemmas.

\begin{lem} \label{special1}
Both surfaces $S$ and $S'$ are normal.
\end{lem}
\dem Since $S$ and $S'$ are hypersurfaces of $K^3$, it suffices to show that they
are both nonsingular in codimension 1 by Serre's criterion. Since $(0,0,0)$ is their only singular point,
they are normal.
\qed

\begin{lem} \label{special2}
The morphism $F$ is the quotient map of the $\Z/2\Z$-action on $S'$
defined by the involution $\sigma$. In particular,
$S=S'/\!/(\Z/2\Z)$.
\end{lem}
\dem By means of the morphism $F$, we can identify ${\cal{O}}(S)$
with a subring of ${\cal{O}}(S')$. Via this identification, we only
need to show that ${\cal{O}}(S)$ is the ring of invariants of
$\sigma^*$. We proceed in several steps. First, we show that $F$ is
a finite morphism. Note that $\alpha^2$ and $\beta +\gamma$ belong
to ${\cal{O}}(S)$, so that $\alpha$ and $\beta + \gamma$ are
integral over ${\cal{O}}(S)$. Since $\alpha^{2e_1} - 4\beta
\gamma=0$, we obtain the relation:
$$
(\beta - \gamma)^2 = - \alpha^{2e_1} + (\beta + \gamma)^2\; .
$$
In particular, $\beta - \gamma$ is integral over ${\cal{O}}(S)$.
Therefore, $\alpha$, $\beta$ and $\gamma$ are integral over
${\cal{O}}(S)$, hence the morphism $F$ is finite. Second, denote by
$K(S)$ (resp. $K(S')$) the field of rational functions of $S$ (resp.
$S'$). We are going to prove that $K(S)$ is the field
$K(S')^{\sigma^*}$ of invariants of $\sigma^*$. It is easy to check
that the functions $\alpha^2$, $\beta + \gamma$ and $\alpha(\beta -
\gamma)$ are $\sigma^*$-invariant. In particular, $K(S)$ is
contained in $K(S')^{\sigma^*}$. Since $\alpha^2$ belongs to $K(S)$
and that $K(S')=K(S)[\alpha]$, $K(S')/K(S)$ is an extension of
degree 2. In particular, we find that $K(S')^{\sigma^*}=K(S)$. Now,
we can prove that  ${\cal{O}}(S')^{\sigma^*}={\cal{O}}(S)$. Since
the coordinate functions of $F$ are $\sigma^*$-invariant,
${\cal{O}}(S)$ is contained in ${\cal{O}}(S')^{\sigma^*}$.
Conversely, let $f$ be any $\sigma^*$-invariant regular function on
$S'$. Since $K(S')^{\sigma^*}=K(S)$, $f$ belongs to $K(S)$. Since
$F$ is finite, $f$ is integral over ${\cal{O}}(S)$. But $S$ is
normal, so ${\cal{O}}(S)$ is integrally closed and $f$ belongs to
${\cal{O}}(S)$. \qed

\begin{lem} \label{special3}
Let $D$ be any derivation on ${\cal{O}}(S)$. Then there exists a unique derivation
$D'$ on ${\cal{O}}(S')$, which commutes with $\sigma^*$ and such that $D'$ coincides
with $D$ on ${\cal{O}}(S)$. Moreover, if $D$ is weighted homogeneous on $S$, then
$D'$ is also weighted homogeneous on $S'$.
\end{lem}
\dem By assumption, $D$ is a $K$-derivation on ${\cal{O}}(S)$, hence it extends uniquely into
a $K$-derivation of $K(S)$. Since $K(S')/K(S)$ is a finite extension, $D$ extends
uniquely into a $K$-derivation $D'$ of $K(S')$ by some elementary results of
Differential Galois Theory (see \cite{Ko}). Consider the map $D_0 = \sigma^* \circ
D' \circ \sigma^*$ on $K(S')$. By construction, $D_0$ is a $K$-linear map on $K(S')$.
Since $\sigma^*$ is an involution, we have for any elements $f$ and $g$ of $K(S')$:
$$
\begin{array}{lll}
D_0(fg) & = & \sigma^*(D(\sigma^*(f)\sigma^*(g))) \\
 & = & \sigma^*(\sigma^*(f)D(\sigma^*(g)) + \sigma^*(\sigma^*(g) D(\sigma^*(f))) \\
 & = & fD_0(g) + gD_0(f)\; .
\end{array}
$$
Therefore, $D_0$ satisfies the Leibniz rule, hence it is a $K$-derivation of $K(S')$. Since $\sigma^*$ is the identity on
$K(S)$, $D$ coincides with $D_0$ on $K(S)$. By uniqueness of the extension of a $K$-derivation
in an algebraic extension, we find that $D=D_0$, hence $D' \circ \sigma^* = \sigma^* \circ D'$. So there only remains
to prove that $D'$ preserves the ring ${\cal{O}}(S')$.

Let $p$ be any point of $S'$ distinct from $(0,0,0)$. Then the point
$q=F(p)$ is distinct from $(0,0,0)$. In particular, $p$ (resp. $q$)
is a smooth point of $S'$ (resp. $S$). Since $F$ is a quotient map
of a finite group of order 2, its singular points correspond to its
fixpoints. Since $(0,0,0)$  is the only fixpoint of $\sigma^*$, $F$
is nonsingular at $p$. In particular, the pull-back morphism $F^*$
defines an isomorphism between the completions of the local rings
${\cal{O}}_{S',p}$ and ${\cal{O}}_{S,q}$ with respect to their
maximal ideals. Via this isomorphism, $D$ determines a unique
derivation $D''$ on the completion of ${\cal{O}}_{S',p}$. By
uniqueness of the extension of $D$, we deduce that $D'=D''$. In
particular, $D'$ is regular at the point $p$. Since this holds for
any point $p\not=(0,0,0)$, $D'$ is regular on $S'-\{(0,0,0)\}$.
Thus, for any regular function $f$ on $S'$, $D'(f)$ is regular on
$S'-\{(0,0,0)\}$. But since $S'$ is normal, $D'(f)$ is regular on
all of $S'$. In particular, $D'$ preserves the ring ${\cal{O}}(S')$.

Finally, assume that $D$ is weighted homogeneous on $S$. By
definition, this means there exists a constant $r$ such that $D$
maps every weighted homogeneous element of degree $n$ to a weighted
homogeneous element of degree $n+r$. Since $\alpha^2$ is a weighted
homogeneous element of ${\cal{O}}(S)$ of degree $2$,
$D'(\alpha)=D(\alpha^2)/2\alpha$ is weighted homogeneous of degree
$r+1=r+\deg(\alpha)$. Since $\beta + \gamma$ is a weighted
homogeneous element of ${\cal{O}}(S)$ of degree $e_1$, $D'(\beta +
\gamma)= D(\beta + \gamma)$ is weighted homogeneous of degree $r +
e_1 =r + \deg(\beta +\gamma)$. Similarly, one can check that
$D'(\beta - \gamma)$ is weighted homogeneous of degree $e_1 + r=r +
\deg(\beta - \gamma)$. Since $\alpha$, $\beta+\gamma$ and $\beta -
\gamma$ span the graded ring ${\cal{O}}(S')$, $D'$ is a weighted
homogeneous $K$-derivation of degree $r$. \qed {\it Proof of
Proposition \ref{special}}: Assume that the polynomial
$R=x_3^2+(x_1^{e_1}-x_2^{2})x_1$ is annihilated  by a nonzero
locally nilpotent derivation. By Lemma \ref{LND1}, the surface $S$
admits a nonzero weighted homogeneous locally nilpotent derivation
$D$.  By Lemma \ref{special3}, $D$ extends uniquely into a weighted
homogeneous $K$-derivation $D'$ of $S'$, which commutes with the
involution $\sigma^*$. Since ${\cal{O}}(S')$ is integral over
${\cal{O}}(S)$ by Lemma \ref{special2} and that $D$ is locally
nilpotent, $D'$ is also locally nilpotent by Vasconcelos' Theorem
(see \cite{Vas}). Set $\varphi(x)=x^{2e_1} /4$. For any polynomial
$f$ of $K[x]$, define as in \cite{Dai} the automorphism $\Delta_f$
of $S'$ by the formula:
$$
\Delta_f(\alpha,\beta,\gamma)= \left (\alpha + \beta f(\beta),
\beta, \frac{\varphi(\alpha +  \beta f(\beta))}{\beta}\right )\; .
$$
Let $\delta$ be the automorphim of $S'$ given by
$\delta(\alpha,\beta,\gamma)=(\alpha,\gamma,\beta)$. Denote by $G$
the group spanned by $\delta$ and all the $\Delta_f$, where $f$ runs
through $K[x]$. Then $G$ acts transitively on the kernels of nonzero
locally nilpotent derivations on $S'$ (see \cite{Dai}). Since
$K[\beta]$ is the kernel of a locally nilpotent derivation on $S'$,
this means that $\ker D'$ must be of the form $K[g]$, where $g$ is
the second coordinate function of an element of $G$. Note that every
generator of $G$ fixes the origin $(0,0,0)$, and that its linear
part at $(0,0,0)$ is either of the form
$(\alpha,\beta,\gamma)\mapsto (\alpha + a \beta +b\gamma,
\beta,\gamma)$ or $(\alpha,\beta,\gamma)\mapsto (\alpha + a \beta
+b\gamma,\gamma,\beta)$. So every element of $G$ fixes the point
$(0,0,0)$, and its linear part at $(0,0,0)$ is either of the form
$(\alpha,\beta,\gamma)\mapsto (\alpha + a \beta +b\gamma,
\beta,\gamma)$ or $(\alpha,\beta,\gamma)\mapsto (\alpha + a \beta
+b\gamma,\gamma,\beta)$. Therefore, $g$ must be either of the form
$\beta +h$ or $\gamma+h$, where $h$ is a polynomial with no constant
nor linear terms in $\alpha,\beta,\gamma$. Since $D'$ is weighted
homogeneous and that $g(0,0,0)=0$, $g$ must be also weighted
homogeneous. So $g$ is either of the form $\beta +\lambda
\alpha^{e_1}$  or $\gamma+\lambda \alpha^{e_1}$, where $\lambda$
belongs to $K$. Now, since $D'$ commutes with $\sigma^*$, $\ker D'$
must be stable by the action of $\sigma^*$. In particular, $g$ is
$\sigma^*$-equivariant. But this is impossible since $\sigma^*$
permutes $\beta$ and $\gamma$ and that $\sigma^*(\alpha)=-\alpha$.
This ends the proof of the proposition.
\qed
\begin{rem}
One has $S'\simeq  K^2/\!/(\Z/2e_1\Z)$ where $\Z/2e_1\Z$ acts on
$K^2$ through $(x,y)\mapsto (-\xi x,-\xi^{-1}y )$, $\xi$ being a
primitive $e_1$-th root of the unity. In  turn, the action on $K^2$
induced by $(x,y)\mapsto (\xi x,\xi^{-1}y)$ and $(x,y)\mapsto
(-y,x)$ gives $S\simeq S'/\!/(\Z/2\Z)=K^2/\!/G$ where $G$ is a
semi-direct product $\Z/e_1\Z\ltimes
 \Z/4\Z$. In case $K=\C$ the fact that $S$ does not admit a
 non-trivial LND ($\C_+$-action) follows from Theorem 2.13 in
 \cite{FZ}: $G$ is not cyclic and the quotient $\C^2/\!/G$ is surely
 not isomorphic to a quotient of $\C^2$ by a cyclic group (the group by which one takes the quotient
 reflects in the fundamental group of the smooth locus of the surface).
\end{rem}


\section{Proof of the third Theorem}\label{idpr}

Assume that the ideal $I$ of relations is nonzero principal, i.e.
$I=(R)$ with $R\not=0$.  By Theorem \ref{relations}, Lemma
\ref{LND1}, Proposition \ref{noLND} and \ref{special} $R$ is
proportional to none of the  polynomials in the following
"forbidden" list:

$$
\begin{array}{ccc}
x_3^2+ax_1^{4}+bx_2^{3} & \mbox{where} & ab\neq 0 \\
x_3^2+ax_1^{5}+bx_2^{3} & \mbox{where} & ab\neq 0\\
 x_3^2+(a_1x_1^{e_1}+b_1x_2)(a_2x_1^{e_1}+b_2x_2)x_1 & \mbox{where} & a_1b_2 - a_2b_1\not=0\; e_1\geq 1\\
x_3^2+(a_1x_1+b_1x_2)(a_2x_1+b_2x_2)(a_3x_1+b_3x_2) & \mbox{where} & a_ib_j - a_jb_i\not=0, \forall i\not=j \\
x_3 ^2 + (ax_1 ^3 + bx_2 ^2)x_2 & \mbox{where} & ab\neq 0\\
x_3 ^2 + (ax_1 ^{e_1} + bx_2 ^2)x_1 & \mbox{where} & ab\neq 0,\; e_1\in 3+2\N\; .
\end{array}
$$

Being homogeneous, the polynomial $R$ can be written as:
$$
R=\sum_{\alpha_1 d_1 + \alpha_2 d_2 + \alpha_3d_3 = \deg_2(R)}
R_{\alpha_1,\alpha_2,\alpha_3} x_1 ^{\alpha_1}  x_2 ^{\alpha_2} x_3
^{\alpha_3}\; .
$$
Denote by $\Supp(R)$ the set
$\{\alpha=(\alpha_1,\alpha_2,\alpha_3)|R_{\alpha}\neq 0\}$. Then
Theorem \ref{inegalite}  implies:
\begin{eqnarray}\label{Louis}
\forall \alpha \in \Supp(R), \quad
\alpha_1d_1+\alpha_2d_2+\alpha_3d_3 \leq d_1+d_2+d_3-2
\end{eqnarray}
Assume that $d_1\leq d_2 \leq d_3$. Since $d_1+d_2+d_3-2<3d_3$, we
find that $\alpha_3\leq 2$ for any element $\alpha$ of $\Supp(
 R)$. In particular, $R$ is a polynomial of degree $\leq 2$ in the
variable $x_3$. Set $d=\gcd(d_1,d_2)$ and write $d_1=de_2$,
$d_2=de_1$. Note that $\gcd(e_1,e_2)=1$ and $e_2\leq e_1$, with
equality if and only if $e_1=e_2=1$.

\subsection{The case when $\deg_{x_3}(R)=0$}

Then, $R$ is an irreducible weighted homogeneous polynomial in the
variables $x_1,x_2$ and for the weights $d_1,d_2$. Since $K$ is
algebraically closed, $R$ is either of the form
$ax_1^{e_1}+bx_2^{e_2}$ or  $cx_i$, where $i\in \{1,2\}$ and
$a,b,c\neq 0$. Since neither $\overline f_1$ nor $\overline f_2$ is
zero, $R$ has type (\ref{t2}).

\subsection{The case when $\deg_{x_3}(R)=1$}

Let $\alpha$ be any point of $\Supp(R)$ of the form
$\alpha=(\alpha_1,\alpha_2,1)$. By inequality (\ref{Louis}), we
have:
$$
\begin{array}{ccl}
\alpha_1d_1+\alpha_2d_2+d_3\leq d_1+d_2+d_3-2 & \Rightarrow &
\alpha_1d_1+\alpha_2d_2\leq d_1+d_2-2 \\
 & \Rightarrow & \alpha_2 \leq 1 \; \mbox{and} \; (\alpha_2=1 \Rightarrow \alpha_1=0)\; .
\end{array}
$$
If $(0,1,1)$ does not belong to $\Supp(R)$, then $\alpha_2$ is
always zero and $R$ has type (\ref{t4}) with $x'_3=x_3$. Now assume
that $(0,1,1)$ belongs to $\Supp(R)$. Then, up to a multiplication
by a constant, $R$ is of the form:
$$
R=(x_2+ax_1^{e_1})x_3+P(x_1,x_2)
$$
where $a\neq 0$ only if $e_2=1$, and where $P(x_1,x_2)$ is weighted
homogeneous. Write $P$ as
$P(x_1,x_2)=P_0(x_1)+(x_2+ax_1^{e_1})Q(x_1,x_2+ax_1^{e_1})$. Taking
$x'_3=x_3+Q(x_1,x_2+ax_1^{e_1})$, we get a polynomial of type
(\ref{t3}).

\subsection{The case when $\deg_{x_3}(R)=2$}

Let $\alpha$ be any element of $\Supp(R)$ of the form
$(\alpha_1,\alpha_2,2)$. Assume that either $\alpha_1$ or $\alpha_2$
is nonzero. Then, by inequality (\ref{Louis}), we would find:
$$
d_1+2d_3\leq\alpha_1d_1+\alpha_2d_2+2d_3\leq d_1+d_2+d_3-2\quad \Rightarrow \quad d_3\leq d_2-2
$$
which is impossible. So, if $\alpha=(\alpha_1,\alpha_2,2)$, then
$\alpha$ is equal to $(0,0,2)$ and
\begin{eqnarray}
2d_3 & \leq & d_1+d_2+d_3-2 \Rightarrow\nonumber \\
d_3 & \leq & d_1+d_2-2  \; .\label{Carla}
\end{eqnarray}

 This means in particular that, up to a scalar multiplication, $R$ is of the
form:
$$
R=x_3^2+x_3P(x_1,x_2)+Q(x_1,x_2)
$$
where $P$ and $Q$ are weighted homogeneous. Taking $x'_3=x_3 +
P(x_1,x_2)/2$, we have:
\begin{eqnarray*}
R & = &
 {x'_3}^2+\mbox{a weighted homogeneous polynomial in }x_1,x_2\\
 & = &
 {x'_3}^2+c\prod_{i=1}^k(a_ix_1^{e_1}+b_ix_2^{e_2})x_1^{r_1}x_2^{r_2}
\end{eqnarray*}
where $r_l<e_l$ for  $l=1,2$ 
 ($c$ is assumed to be one as soon as $k\geq 1$). Since $R$
is weighted homogeneous, we obtain:
\begin{eqnarray*}
2d_3 & = & ke_kd_k+r_1d_1+r_2d_2 \mbox{ for }k=1,2\\
     & = & kde_1e_2+r_1de_2+r_2de_1\; .
\end{eqnarray*}
Inequality (\ref{Carla}) yields $d_3\leq de_1+de_2-2$ and since
$d_2\leq d_3$ one gets:
$$
2d_2= 2de_1\leq 2d_3=kde_1e_2+r_1de_2+r_2de_1\leq 2de_1+2de_2-4\; .
$$
After division by $de_1e_2$, this yields:
\begin{eqnarray}
2/e_2\leq k+r_1/e_1+r_2/e_2<2/e_2+2/e_1\; .\label{Marine}
\end{eqnarray}
In particular, we have $k<4$. We are going to describe the different
possible polynomials $R$ depending on the value of $k$. \\
{\underline{$1^{st}$ case: $k=0$} i.e. $R={x'_3}^2+cx_1^{r_1}x_2^{r_2}$ . \\
Inequalities (\ref{Marine}) becomes $2/e_2\leq
r_1/e_1+r_2/e_2<2/e_2+2/e_1$ and we deduce that either $r_1\leq 1$
or $r_2\leq 1$. \\
If $r_1=0$, then we have:
$$
2/e_2\leq r_2/e_2<2/e_2+2/e_1\; .
$$
After multiplication by $e_2$, we find $2\leq r_2<2+2e_2/e_1$. Since
$e_2/e_1\leq 1$, we obtain that either $r_2=2$  or $r_2=3$. Since
$R$ is irreducible, the case $r_2=2$ is impossible and we have:
$$
R={x'_3}^2+cx_2^{3}\; .
$$
So $R$ corresponds to a polynomial of type (\ref{t5}) in the list.\\
If $r_1=1$, then we have:
$$
\begin{array}{ccccccc}
2/e_2\leq 1/e_1+r_2/e_2<2/e_2+2/e_1 & \Rightarrow  & 2e_1/e_2 & \leq &  1+r_2e_1/e_2 & <& 2+2e_1/e_2 \\
& & \\
& \Rightarrow & 2e_1/e_2-1 & \leq &   r_2e_1/e_2  & < &2e_1/e_2+1 \\
& & \\
& \Rightarrow & -1 & \leq & (r_2-2)e_1/e_2 & < &1\; .
\end{array}
$$
Since  $e_1/e_2\geq 1$, $r_2$ is either equal to $1$ or $2$.
Therefore $R$ is either of type (\ref{t8}) or (\ref{t6}).\\
 Now, if
$r_1\geq 2$, then $r_2$ is either equal to $0$ or $1$. If $r_2=0$,
then, by irreducibility, $r_1\in 3+2\N$ and $R$ is of type
(\ref{t7}).
If $r_2=1$, then $R$ is of type (\ref{t8}). \\ \\
{\underline{$2^{nd}$ case: $k=1$}} i.e. $R={x'_3}^2+(a{x_1}^{e_1}+b{x_2}^{e_2})x_1^{r_1}x_2^{r_2}$.\\
Inequalities (\ref{Marine}) become $2/e_2\leq
1+r_1/e_1+r_2/e_2<2/e_2+2/e_1$ and it follows that either $r_1\leq
1$ or $r_2\leq 1$.\\
 If $r_1=0$, then we have:
$$
2/e_2\leq 1+r_2/e_2<2/e_2+2/e_1\; .
$$
After multiplication by $e_2$, we find $2\leq e_2+r_2<2+2e_2/e_1$.
Since $e_2/e_1 \leq 1$, $e_2+r_2$ is either equal  to $2$ or $3$.
Assume that $e_2+r_2=2$. Since $r_2<e_2$, we have $r_2=0$ and
$e_2=2$. Therefore, $R$ is of type (\ref{t9}). Assume now that
$e_2+r_2=3$. Then we find that $1<2e_2/e_1$, hence $e_1<2e_2$.
Recalling that $e_2\leq e_1$ with equality only if $e_1=e_2=1$, we are left with the following three
cases:
$$
\left \{
\begin{array}{cc}
k=1,\,r_1=0,\,r_2=1,e_2=2,e_1=3 &\\
k=1,\,r_1=0,\,r_2=0,e_2=3,e_1=4 &\\
k=1,\,r_1=0,\,r_2=0,e_2=3,e_1=5 &\; .
\end{array}
\right.
$$
Therefore, $R$ is equal to one of the following polynomials:
$$
\left \{
\begin{array}{lc}
{x'_3}^2+(ax_1^3+bx_2^2)x_2&\\
{x'_3}^2+ax_1^4+bx_2^3&\\
{x'_3}^2+ax_1^5+bx_2^3\; .&
\end{array}
\right.
$$
They all belong to the forbidden list at the beginning of this
section unless $ab=0$. Since ${x'_3}^2+ax_1^4$ is reducible, $R$ has
either type (\ref{t8}),
(\ref{t5}) or (\ref{t7}).\\
 If now $r_1=1$, then we have:
$$
2/e_2\leq 1+1/e_1+r_2/e_2<2/e_2+2/e_1\; .
$$
After multiplication by $e_2$, we obtain:
$$
2\leq e_2+r_2+\frac{e_2}{e_1}<2+2\frac{e_2}{e_1} \quad \Rightarrow \quad
2-\frac{e_2}{e_1}\leq e_2+r_2<2+\frac{e_2}{e_1}\; .
$$
Since $e_2/e_1\leq 1$, we find that $e_2+r_2$ is either equal to $1$
or $2$. Assume that $e_2+r_2=1$.  Then $2 - e_2/e_1=1$ i.e. $e_1=e_2=1$ and, since
$r_2<e_2$, we obtain that $r_2=0$. Therefore, $R$ is
of the form ${x'_3}^2+(ax_1+bx_2)x_1$, i.e. it is of type
(\ref{t8.5}). Assume now that $e_2 + r_2=2$. Since $r_2 <e_2$, we
find that $r_2=0$ and $e_2=2$. Therefore, recalling that $k=1$ and $e_1>r_1=1$, $R$ is
of the form $ R={x'_3}^2+(ax_1^{e_1}+bx_2^2)x_1$ with
$1=\gcd(e_1,e_2)=\gcd(e_1,2)\Rightarrow e_1\in 3+2\N$. So $R$
belongs to the forbidden list unless $ab=0$. If $a=0$ and $b\neq 0$,
then $R$ is of type (\ref{t6}). If $a\neq0$ and $b=0$, then $R$ is
reducible.\\
Finally assume $r_1\geq 2$. The right hand side inequality in
(\ref{Marine}) yields
$$
\begin{array}{ccccccc}
1+2/e_1+r_2/e_2 & < & 2/e_2+2/e_1 & \Rightarrow & 1+r_2/e_2 & < & 2/e_2 \\
 & & & \Rightarrow & e_2+r_2 & < & 2\; .
\end{array}
$$
Since $r_2<e_2$, we get that $k=1$, $r_1\geq 2$, $r_2=0$ and
$e_2=1$. Therefore,
$R$ is of the form ${x'_3}^2+(ax_1^{e_1}+bx_2)x_1^{r_1}$, i.e. it is of type (\ref{t8.5}). \\ \\
{\underline{$3^{rd}$ case: $k=2$} i.e.
$R={x'_3}^2+(a_1{x_1}^{e_1}+b_1{x_2}^{e_2})(a_2{x_1}^{e_1}
+b_2{x_2}^{e_2})x_1^{r_1}x_2^{r_2}$.\\
 We start from (\ref{Marine}):
$$
\begin{array}{cccccccclcc}
2/e_2 & \leq & 2+r_1/e_1+r_2/e_2 & < & 2/e_2+2/e_1 & \Rightarrow & 2 & < & 2/e_2+2/e_1  & \leq & 4/e_2 \\
& & & & & \Rightarrow & e_2 & < & 2 & & \\
\end{array}
$$
so we get $ e_2=1$, $r_2=0$ and, using again the left hand side of
the implication: $2\leq 2+r_1/e_1<2+2/e_1$ which simply means
$r_1<2$. Therefore, $R$ is equal to one of the following
polynomials:
$$
\left \{
\begin{array}{lc}
{x'_3}^2+(a_1x_1^{e_1}+b_1x_2)(a_2x_1^{e_1}+b_2x_2) & \\
{x'_3}^2+(a_1x_1^{e_1}+b_1x_2)(a_2x_1^{e_1}+b_2x_2)x_1 &\; .
\end{array}
\right.
$$
The first one is of type (\ref{t3.5}) when not reducible.\\
The second one is in the forbidden list unless $a_1b_2-a_2b_1=0$,
which gives a polynomial of the form ${x'_3}^2+(ax_1^{e_1}+bx_2)^2
x_1$, i.e. of type
(\ref{t14}).  \\
{\underline{$4^{th}$ case: $k=3$}  i.e.
$R={x'_3}^2+(a_1{x_1}^{e_1}+b_1{x_2}^{e_2})(a_2{x_1}^{e_1}
+b_2{x_2}^{e_2})(a_3{x_1}^{e_1}+b_3{x_2}^{e_2})x_1^{r_1}x_2^{r_2}$. \\
Again (\ref{Marine}) yields $$ 2/e_2 \leq  3+r_1/e_1+r_2/e_2  <
2/e_2+2/e_1 \quad \Rightarrow \quad 3  <  2/e_2+2/e_1
$$
which is possible only if  $e_1=e_2=1$. Consequently $r_1=r_2=0$ and
$R$ is of the form:
$$
R= {x'_3}^2+(a_1x_1+b_1x_2)(a_2x_1+b_2x_2)(a_3x_1+b_3x_2)\; .
$$
This is again in the forbidden list unless one of the determinants
$a_ib_j - a_jb_i$ is zero. In this case, $R$ is of type
(\ref{t6.5}), and the classification is complete.

\section{Proof of the Corollary }\label{democor}
{\it Proof of Corollary \ref{THEcor}}:  By Theorem \ref{fleche} we
have to prove that every one of the 14 polynomials in the list given
inside this theorem is tamely equivalent to $0$, $x_3$,
$x_1^r+x_2^s$ or $x_1^kx_3+P(x_1,x_2)$.
Let us introduce the notation $Q\equiv R$  when $Q$ and $R$ are tamely equivalent and we call $R_i$,
  $\forall i=1,\cdots,14$, the polynomial such the relation (i) in the list is equal to $R_i(x_1,x_2,x_3')$.
  Applying the elementary automorphism $(x_1,x_2,x_3-h(x_1,x_2))$ we get $R_i(x_1,x_2,x_3')\equiv R_i(x_1,x_2,x_3)$
  so we may forget the "prime" in $x_3'$.  \\
We have $R_1=0$ and $R_2=x_3$ so there is  nothing to do. One has
$R_3\circ(x_1,dx_2,x_3)=x_1^{e_1}+x_2^{e_2}$ where $d^s=1/c$ so
$R_3\equiv x_1^r+x_2^s$ with $(r,s)=(e_1,e_2)$. One has
$R_4\circ(x_1,x_2+ax_1^{e_1},x_3)=R_3\equiv x_1^r+x_2^s$ so
$R_4\equiv x_1^r+x_2^s$. The relation $R_5$ is even equal to
$x_1^kx_3+P(x_1,x_2)$. By a transposition $(x_3,x_2,x_1)$ we get
$R_6\equiv R_3\equiv x_1^r+x_2^s$ with $(r,s)=(2,3)$. Similarly,
using if needed some permutation and some tame automorphism of the
form $(x_1,\frac{1}{b}(x_2-a x_1^{e_1}))$ one gets that $R_7$,
$R_9$, $R_{11}$ and $R_{14}$ are tamely equivalent to
$x_1^kx_3+P(x_1,x_2)$ and $R_8$ to $x_1^r+x_2^s$. As for $R_{10}$,
let $d^2=-b$ and compose $R_{10}\circ
(x_1,\frac{x_2-x_3}{2d},\frac{x_2+x_3}{2})=ax_1^{e_1}+x_3x_2\equiv
x_1^kx_3+P(x_1,x_2)$. In $R_{12}$, if $b_1b_2=0$ then clearly
$R_{12}\equiv R_{11}$ or direcly $R_9$, that is, equivalent to
$x_1^kx_3+P(x_1,x_2)$. If $b_1b_2\neq 0$ then one checks that
$R_{12}\circ(x_1,x_2-\frac{a_1b_2+b_1a_2}{2b_1b_2}x_1^{e_1},x_3)$
looks like $R_{10}$ but with an even $e_1$, which does not prevent
us from showing, as before, that $R_{12}\equiv x_1^kx_3+P(x_1,x_2)$.
Finally, using an appropriate linear automorphism shows that
$R_{13}\equiv R_7\equiv x_1^kx_3+P(x_1,x_2)$ unless
$a_1b_2-a_2b_1=0$, in which case $R_{13}\equiv R_3$. \qed Remark
that the four polynomials in Corollary \ref{THEcor} do belong to
kernels of nontrivial LNDs of $K[x_1,x_2,x_3]$: 0 and $x_3$ are in
the kernel of $\frac{\partial}{\partial x_1}$, $x_1^r+x_2^s$ in the
one of $\frac{\partial}{\partial x_3}$ and $x_1^kx_3+P(x_1,x_2)$ is
in the kernel of $x_1^k\frac{\partial}{\partial x_2}-\frac{\partial
P}{\partial x_2}\frac{\partial}{\partial x_3}$.

\end{document}